\documentclass[11pt]{amsart}
\usepackage{amssymb}

\newtheorem{theorem}{Theorem}[section]
\newtheorem{lemma}[theorem]{Lemma}
\theoremstyle{definition}
\newtheorem{definition}[theorem]{Definition}
\newtheorem{prop}[theorem]{Proposition}
\newtheorem{corollary}[theorem]{Corollary}
\theoremstyle{remark}
\newtheorem{remark}[theorem]{Remark}
\numberwithin{equation}{section}

\def\qed{\hfill {\hbox{${\vcenter{\vbox{               %HOLLOW SQUARE
   \hrule height 0.4pt\hbox{\vrule width 0.4pt height 6pt
   \kern5pt\vrule width 0.4pt}\hrule height 0.4pt}}}$}}}
\def\domedskip{\par \ifdim\lastskip<\medskipamount
  \removelastskip \medskip\fi}
\newcommand{\db}[1]{\ensuremath{\bar{\bar{\mbox{\ensuremath{ #1}}}}}}
\begin{document}

%%%%%%%%%%%%%%%%%%%%%%%%%%%%%%%%%%%%%%%%%%%%%%%%%%%%%%%%%%%%
%%%%%%%%%%%%%%%%%%%%%%%%%%%%%%%%%%%%%%%%%%%%%%%%%%%%%%%%%%%%
% This a placeholder for the TOPOLOGY PROCEEDINGS logo %%%%%%
\noindent                                             %%%%%%
\begin{picture}(150,36)                               %%%%%%
\put(5,20){\tiny{Submitted to}}                       %%%%%%
\put(5,7){\textbf{Topology Proceedings}}              %%%%%%
\put(0,0){\framebox(140,34)}                          %%%%%%
\put(2,2){\framebox(136,30)}                          %%%%%%
\end{picture}                                         %%%%%%
%%%%%%%%%%%%%%%%%%%%%%%%%%%%%%%%%%%%%%%%%%%%%%%%%%%%%%%%%%%%
%%%%%%%%%%%%%%%%%%%%%%%%%%%%%%%%%%%%%%%%%%%%%%%%%%%%%%%%%%%%

\vspace{0.5in}

\title[Classification of Finite Alexander Quandles]%
{Classification of Finite Alexander Quandles}

\author{Sam Nelson}
\address{Department of Mathematics, Louisiana State University,
Baton Rouge, LA 70803}
\email{nelson@math.lsu.edu}



\begin{abstract}
Two finite Alexander quandles with the same number of elements are 
isomorphic iff their $\mathbb{Z}[t^{\pm 1}]$-submodules $\mathrm{Im}(1-t)$ 
are isomorphic as modules. This yields specific conditions on 
when Alexander quandles of the form $\mathbb{Z}_n[t^{\pm 1}]/(t-a)$ where
$\gcd(n,a)=1$ (called \textit{linear} quandles) are isomorphic, as well 
as specific conditions on when two linear quandles are dual and which linear 
quandles are connected. We apply this result, obtaining a procedure for 
classifying Alexander quandles of any finite order and as an application 
we list the numbers of distinct and connected Alexander quandles with up 
to fifteen elements.
\end{abstract}

\maketitle

\section{Introduction}

In \cite{J}, D. Joyce defined the fundamental quandle, an algebraic 
invariant of knots which classifies classical knots. The set of quandles 
forms a category whose axioms are algebraic versions of the three 
Reidemeister moves. Quandles are useful both for defining new knot 
invariants (as in \cite{CJKS}) and for improving our understanding 
of old ones (see \cite{FR}, for example). 

The ability to distinguish quandles would allow us to distinguish knots.
While there is not yet a complete classification theorem for general 
quandles, there are classification results for quandles of prime order
\cite{O} and for indecomposable quandles of prime squared order \cite{G}.
In this paper we classify finite Alexander quandles by reducing the problem
of comparing finite Alexander quandles to comparing certain 
$\mathbb{Z}[t^{\pm 1}]$-submodules.

\begin{definition}
A \textit{quandle} is a set $X$ with a binary operation written as 
exponentiation satisfying 

\begin{list}{}{}
\item[(i)] For every $a,b \in X$ there exists a unique $c\in X$ such 
that $a=c^b$,
\item[(ii)] For every $a,b,c\in X$ we have $a^{bc}=a^{cb^c}$, and
\item[(iii)] For every $a\in X$ we have $a^a=a$.
\end{list}
\end{definition}

Any module over $\Lambda = \mathbb{Z}[t^{\pm 1}]$ is a quandle under 
the operation $a^b=ta+(1-t)b$. Quandles of this form are called 
\textit{Alexander quandles}. To obtain finite Alexander quandles, we 
typically consider $\Lambda_n /(h)$ where 
$\Lambda_n=\mathbb{Z}_n[t^{\pm 1}]$ and $h$ is a monic polynomial in 
$t$. In an earlier version of \cite{O}, the questions of
when two Alexander quandles of the form $\Lambda_n/(t-a)$ with $\gcd(n,a)=1$ 
(we call Alexander quandles of this form \textit{linear}) are 
isomorphic and of when two linear quandles are dual were posed.

To answer these questions, we first consider the general case of when
two arbitrary Alexander quandles of finite cardinality are isomorphic. 
We obtain a result which reduces the problem of comparing Alexander 
quandles to comparing certain $\Lambda$-submodules. We then apply this 
result, obtaining a pair of simple conditions on $a$ and $b$ which are 
necessary and sufficient for two linear Alexander quandles 
$\Lambda_n/(t-a)$ and $\Lambda_n/(t-b)$ to be isomorphic.

In the course of answering the question of classifying linear quandles, 
we also answer the question of when linear quandles are dual and we
obtain results on when Alexander quandles are connected.

\section{Alexander quandles and $\Lambda$-modules}

Since the quandle structure of an Alexander quandle is determined by
its $\Lambda$-module structure, any isomorphism of $\Lambda$-modules 
is also an isomorphism of Alexander quandles. The converse is not true, 
however: $\Lambda_9/(t-4)$ is isomorphic to $\Lambda_9/(t-7)$ as an
Alexander quandle but not as a $\Lambda$-module.

Nonetheless, an isomorphism of Alexander quandles is in a sense almost
an isomorphism of $\Lambda$-modules; in fact, (after applying a shift
if necessary) the restriction of a quandle isomorphism $f:M\to N$ to the 
submodule $(1-t)M$ is a $\Lambda$-module isomorphism onto the image
of the restriction. Theorem \ref{thm:main} says that the converse is true 
as well; that is, we can determine whether two Alexander quandles of the 
same finite cardinality are isomorphic simply by comparing these 
$\Lambda$-submodules. This reduces the problem of classifying finite 
Alexander quandles to comparing $\Lambda$-modules of the form $(1-t)M$.

\begin{theorem} \label{thm:main}
Two finite Alexander quandles $M$ and $N$ of the same cardinality are
isomorphic as quandles iff there is an isomorphism of $\Lambda$-modules 
$h:(1-t)M\to(1-t)N$.
\end{theorem}

\begin{proof}
Let $M$ and $N$ be finite Alexander quandles and $f:M\to N$ a quandle 
isomorphism. We may assume 
without loss of generality that $f(0)=0$ since $f':M\to N$ defined by 
$f'(x)=f(x)+c$ is also an isomorphism of Alexander quandles for any
$c\in N$. Then $f(tx+(1-t)y) = tf(x)+(1-t)f(y)$ implies
\[
f(tx)=f(tx+(1-t)0)=tf(x)+(1-t)f(0)=tf(x)
\]
and
\[
f((1-t)y)=f(t0+(1-t)y)=tf(0)+(1-t)f(y)=(1-t)f(y)
\]
so that
\begin{equation}
f(tx+(1-t)y)=f(tx)+f((1-t)y).
\end{equation}

Denote $M'=(1-t)M$ and $N'=(1-t)N$. Since $t^{-1}\in\Lambda$, every 
element of $M$ is $tx$ for some $x\in M$, and since $f(0)=0$, $f$ takes
the coset $0+M'$ of $M'$ in $\bar{M}=M/M'$ to the coset $0+N'$ of $N'$ 
in $\bar{N}=N/N'$, so we have that $h=f|_{M'}:M'\to N'$ is a 
homomorphism of $\Lambda$-modules. Since $f$ is injective, its 
restriction $h$ is a bijection onto its image $0+N'=N'$, and hence $h$ is
an isomorphism of $\Lambda$-modules.

Conversely, suppose $h:M' \to N'$ is an isomorphism of finite 
$\Lambda$-modules with $|M|=|N|$. Let $A\subset M$ be a set of 
representatives of cosets of $M'$ in $\bar{M}$. Then every $m\in M$ has 
the form $m=\alpha + \omega$ for a unique $\alpha \in A$ and 
$\omega \in M'$. We will show that there exists a bijection $k:A \to B$
onto a set $B$ of representatives of cosets of $N'$ in $\bar{N}$ such that
the map $f:M\to N$ defined by 
\[
f(\alpha + \omega) = k(\alpha) + h(\omega)
\]
is an isomorphism of Alexander quandles (though typically not of 
$\Lambda$-modules).

Let $\alpha_1,\alpha_2\in A$ and $\omega_1,\omega_2\in (1-t)M$.
For any $\alpha_1\in A$ we have $t\alpha_1=\alpha_1-
(1-t)\alpha_1$, so that
\begin{eqnarray*}
 &   & \hskip -0.75in f(t(\alpha_1+\omega_1)+(1-t)(\alpha_2+\omega_2)) \\
 & = & f(\alpha_1+t\omega_1+(1-t)(\alpha_2-\alpha_1+\omega_2)) \\
 & = & k(\alpha_1) + h(t\omega_1+(1-t)(\alpha_2-\alpha_1+\omega_2)) \\
 & = & k(\alpha_1) + th(\omega_1)+h((1-t)\alpha_2) \\
 &   & \quad{} -h((1-t)\alpha_1)+(1-t)h(\omega_2). \\
\end{eqnarray*}
On the other hand, 
\begin{eqnarray*}
 &   & \hskip -0.5in tf(\alpha_1+\omega_1)+(1-t)f(\alpha_2+\omega_2) \\
 & = & t(k(\alpha_1)+h(\omega_1)) + (1-t)(k(\alpha_2)+h(\omega_2)) \\
 & = & tk(\alpha_1)+th(\omega_1)+(1-t)k(\alpha_2)+(1-t)h(\omega_2)
\end{eqnarray*}
so for $f$ to be a homomorphism of quandles it is sufficient that
\begin{equation} \label{eqn:nec}
(1-t)k(\alpha_1)-h((1-t)\alpha_1) = (1-t)k(\alpha_2)-h((1-t)\alpha_2)
\end{equation}
for all $\alpha_1,\alpha_2\in A$. We will show that given a set of coset 
representatives $A\subset M$ we can choose a set $B\subset N$ of coset 
representatives and a bijection $k:A\to B$ so that $(1-t)k(\alpha)=
h((1-t)\alpha)$ for all $\alpha\in A$, which satisfies (\ref{eqn:nec}) 
and thus yields a homomorphism $f:M\to N$ of Alexander quandles.
Since this $f$ is setwise the Cartesian product $k\times h$ of the 
bijections $k:A\to B$ and $h:M' \to N'$, $f$ is bijective and hence
an isomorphism of quandles.

Denote $M''=(1-t)^2M$, $\db{M}=M'/M''$ and similarly for $N$. The 
isomorphism $h:M'\to N'$ induces an isomorphism $\bar{h}:\db{M}\to \db{N}$.
There are surjective maps $\psi:\bar{M}\to\db{M}$ and $\phi:\bar{N}\to \db{N}$
induced by multiplication by $(1-t)$. Then $|M'|=|N'|$ and $|M|=|N|$ imply
that $|\bar{M}|=|\bar{N}|$, and in turn $|\db{M}|=|\db{N}|$. 

Then $|\psi^{-1}(y)|=|\psi^{-1}(y')|$ for all $y,y'\in \db{M}$, since 
$\psi(y)=\psi(0)+\psi(y)$ implies that for each element of $\psi^{-1}(0)$ 
there is an element of $\psi^{-1}(y)$, that is, $|\psi^{-1}(y)| \ge 
|\psi^{-1}(0)|$, and similarly $\psi(0)=\psi(0)+\psi(y)+\psi(-y)$ implies that
$|\psi^{-1}(0)| \ge |\psi^{-1}(y)|$. Hence $|\psi^{-1}(y)|=|\bar{M}|/|\db{M}|$
for all $y\in \db{M}$, and similarly $|\phi^{-1}(h(y))|=|\bar{N}|/|\db{N}|=
|\bar{M}|/|\db{M}|$ for all $y\in \db{M}$. Thus there is 
a bijection of sets $g:\bar{M}\to\bar{N}$ such that the diagram
\[
\begin{array}{rcccl}
     & \bar{M}    & \stackrel{g}{\longrightarrow} & \bar{N} & \\
\psi & \downarrow &                             & \downarrow & \phi \\
     & \db{M}    & \stackrel{\bar{h}}{\longrightarrow} & \db{N} & 
\end{array}
\]
commutes.

Let $B$ be a set of coset representatives for $\bar{N}$. Then there is 
a unique bijection $k:A\to B$ such that
\[
\begin{array}{rcccl}
     & A          &\stackrel{k}{\longrightarrow} & B          &      \\
     & \downarrow &                              & \downarrow &      \\
     & \bar{M}    &\stackrel{g}{\longrightarrow} & \bar{N}    &      \\
\psi & \downarrow &                              & \downarrow & \phi \\
     & \db{M}     & \stackrel{\bar{h}}{\longrightarrow} & \db{N} & 
\end{array}
\]
commutes. In particular, $\bar{h}\psi(\alpha) = \phi k(\alpha)$, that is,
\begin{equation} \label{eqn:*}
\bar{h}((1-t)\alpha +(1-t)^2M) =  (1-t)k(\alpha)+(1-t)^2N.
\end{equation}

Define $\gamma:M'\to \db{M}$ and $\epsilon:N'\to \db{N}$ by 
$\gamma((1-t)m)=(1-t)m+(1-t)^2M\in\bar{M}$ and 
$\epsilon((1-t)n)=(1-t)n+(1-t)^2N\in\bar{N}$, the classes of
$(1-t)m$ and $(1-t)n$ in $\db{M}$ and $\db{N}$ respectively. We then
have commutative diagrams
\[
\begin{array}{rcccl}
 & A          & \stackrel{(1-t)}{\longrightarrow} & M'          &           \\
 & \downarrow &                                   & \downarrow  & \! \gamma \\
 & \bar{M}    & \stackrel{\psi}{\longrightarrow}  & \db{M}      &  
\end{array}
\ \ \ \mathrm{and} \ \ \
\begin{array}{rcccl}
 & B          & \stackrel{(1-t)}{\longrightarrow} & N'          &           \\
 & \downarrow &                                 & \downarrow  & \! \epsilon \\
 & \bar{N}    & \stackrel{\phi}{\longrightarrow}  & \db{N}.      &  
\end{array}
\]
Equation (\ref{eqn:*}) then says that outside rectangle of the
diagram
\[
\begin{array}{rcccl}
       & A          &\stackrel{k}{\longrightarrow} & B          &          \\
(1-t)  & \downarrow &                              & \downarrow & (1-t)    \\
       & M'         &\stackrel{h}{\longrightarrow} & N'         &          \\
\gamma & \downarrow &                              & \downarrow & \epsilon \\
       & \db{M}     & \stackrel{\bar{h}}{\longrightarrow} & \db{N} & 
\end{array}
\]
commutes. The bottom square commutes by definition of $\bar{h}$, and thus 
we have $\epsilon(h((1-t)\alpha))=\epsilon((1-t)k(\alpha))$, that is,
\[
h((1-t)\alpha)+(1-t)^2N=(1-t)k(\alpha)+(1-t)^2N.
\]
In particular, there is a $\xi\in N$ so that
\[
h((1-t)\alpha)=(1-t)k(\alpha)+(1-t)^2\xi=(1-t)(k(\alpha)+(1-t)\xi).
\]
Then for each $\alpha\in A$ with $\xi\ne 0$ we may replace $k(\alpha)$ 
with the coset representative $k'(\alpha)=k(\alpha)+(1-t)\xi$ to obtain
a new set $B'$ of coset representatives for $\bar{N}$ and a bijection
$k':A\to B'$ with $(1-t)k'(\alpha)=h((1-t)\alpha)$ so that (\ref{eqn:nec})
is satisfied. Then $f:M\to N$ by $f(\alpha+\omega)=k'(\alpha)+h(\omega)$ 
for all $\alpha\in A$ is an isomorphism of Alexander quandles, as required.
\end{proof}

As a consequence, we obtain Corollary \ref{cor:cond}, which gives specific
conditions on $a$ and $b$ for $\Lambda_n/(t-a)$ and $\Lambda_n/(t-b)$ to be
isomorphic Alexander quandles when $a$ and $b$ are coprime to $n$. 

Denote $N(n,a)=\frac{n}{\gcd(n,1-a)}$ for any $a\in \mathbb{Z}_n$. We
will use the symbol $\cong$ to denote an isomorphism of quandles and
$\approx$ to denote an isomorphism of $\Lambda$-modules.

\begin{corollary} \label{cor:cond}
Let $a$ and $b$ be coprime to $n$. Then the Alexander quandles 
$\Lambda_n/(t-a)$ and $\Lambda_n/(t-b)$ are isomorphic iff
$N(n,a)=N(n,b)$ and $a\equiv b (\mathrm{mod}\ N(a,b))$.
\end{corollary}

\begin{proof}
By theorem \ref{thm:main}, 
\[
\Lambda_n/(t-a)\cong\Lambda_n/(t-b) \iff 
(1-t)[\Lambda_n/(t-a)]\approx(1-t)[\Lambda_n/(t-b)].
\]

As a $\mathbb{Z}$-module, $(1-t)[\Lambda_n/(t-a)]$ is $(1-a)\mathbb{Z}_n$
and $(1-t)[\Lambda_n/(t-b)]$ is $(1-b)\mathbb{Z}_n$ with 
the action of $t$ given by multiplication by $a$ in $(1-a)\mathbb{Z}_n$ 
and by $b$ in $(1-b)\mathbb{Z}_n$.

The $\mathbb{Z}$-module $(1-a)\mathbb{Z}_n$ is isomorphic to 
$\mathbb{Z}_n/\mathrm{Ann}(1-a)$, so 
\begin{eqnarray*}
\Lambda_n/(t-a)\cong\Lambda_n/(t-b) & \iff & 
 \mathbb{Z}_n/\mathrm{Ann}(1-a)\approx \mathbb{Z}_n/\mathrm{Ann}(1-b) \\
 & \iff & \mathrm{Ann}(1-a) = \mathrm{Ann}(1-b) \\
 & \iff & \mathrm{Ord}_{\mathbb{Z}_n}(1-a) =\mathrm{Ord}_{\mathbb{Z}_n}(1-b) \\
 & \iff & \frac{n}{\gcd(n,1-a)}=\frac{n}{\gcd(n,1-b)} \\
 & \iff & N(n,a)=N(n,b).
\end{eqnarray*}
Denote $n'=N(n,a)=N(n,b)$. Then $(1-t)[\Lambda_n/(t-a)]$ is $\mathbb{Z}_{n'}$ 
with $t$ acting by multiplication by $a$, and if $N(n,a)=N(n,b)=n'$
then $(1-t)[\Lambda_n/(t-b)]$ is $\mathbb{Z}_{n'}$ with $t$ acting by 
multiplication by $b$.

Multiplication by $a$ agrees with multiplication by $b$ on $\mathbb{Z}_{n'}$
iff $a\equiv b (\mathrm{mod} \ n')$, so the $\Lambda$-module structures on
$\mathbb{Z}_{n'}$ determined by $a$ and $b$ agree iff 
$a\equiv b (\mathrm{mod} \ n')$.
\end{proof}

\medskip

\begin{definition}
A quandle $M$ is \textit{connected} if it has only one orbit, i.e. if 
the set $\{a^b:b\in M\}=M$ for all $a\in M$. In particular, an 
Alexander quandle is connected if $(1-t)M=M$.
\end{definition}

\begin{corollary} \label{cor:conn}
Two finite connected Alexander quandles are isomorphic iff they are isomorphic
as $\Lambda$-modules.
\end{corollary}

\begin{proof}
This follows from the proof of theorem \ref{thm:main}. Specifically, if
$M$ and $N$ are connected and $f:M\to N$ is an isomorphism of 
quandles with $f(0)=0$, then $f$ is an isomorphism of $\Lambda$-modules.
\end{proof}

\begin{corollary} \label{cor:con2}
A linear Alexander quandle $\Lambda_n/(t-a)$ is connected iff $\gcd(n,1-a)=1$.
\end{corollary}

\begin{proof}
An Alexander quandle is connected iff $M=(1-t)M$. Since 
$(1-t)[\Lambda_n/(t-a)]$ is $\mathbb{Z}_{n_a}$ with t acting 
by multiplication by $a$, we have $\Lambda_n/(t-a)$ is connected iff $n_a=n$, 
that is, iff $\gcd(n,1-a)=1$.
\end{proof}

\begin{corollary} \label{cor:con3}
No linear Alexander quandle $\Lambda_n/(t-a)$ with $n$ even is connected.
\end{corollary}

\begin{proof}
To have a linear quandle $\Lambda_n/(t-a)$ with $n$ elements, we must have
$\gcd(n,a)=1$, so if $n$ is even, $a$ must be odd. But then $1-a$ is even and 
$\gcd(n,1-a)\ne 1$, and $\Lambda_n/(t-a)$ is not connected.
\end{proof}

\medskip

For each $y\in X$ we can define a map of sets  $f_y:X\to X$ by $f_y(x)=x^y$. 
Quandle axiom (i) then says that $f_y$ is a bijection for each
$y\in X$. We may then define a new quandle structure on $X$ by 
$x^{\bar{y}}=f_y^{-1}(x)$; this is the \textit{dual} quandle of $X$.

\begin{lemma} \label{lem:dual}
The dual of an Alexander quandle $X$ is the set $X$ with quandle 
operation given by $x^{\bar{y}}=t^{-1}x+(1-t^{-1})y$. 
\end{lemma}

\begin{proof}
If $f_y(x)=c=tx+(1-t)y$ then $t^{-1}c=x+(t^{-1}-1)y \Rightarrow 
x=t^{-1}c+(1-t^{-1})y$; thus $f_y^{-1}(x)=t^{-1}x+(1-t^{-1})y$.
\end{proof}

\begin{corollary} \label{cor:dual}
Let $a,b$ be coprime to $n$. Then $\Lambda_n/(t-a)$ is dual to 
$\Lambda_n/(t-b)$ iff $n_a=n_b$ and $ab\equiv 1 (\mathrm{mod} \ n_a)$. 
In particular, a linear Alexander quandle $\Lambda_n/(t-a)$ is self-dual 
iff $a$ is a square mod $n_a$.
\end{corollary}

\begin{proof}
If $n$ and $a$ are coprime, then $a$ is invertible in $\mathbb{Z}_n$
and the dual of $\Lambda_n/(t-a)$ is given by $\Lambda_n/(t-a^{-1})$ by
lemma \ref{lem:dual}. Then corollary \ref{cor:cond} says that
$\Lambda_n/(t-b)$ is isomorphic to $\Lambda_n/(t-a^{-1})$ iff 
$n_b=n_{a^{-1}}$ and $b\equiv a^{-1} (\mathrm{mod} \ n_b)$. 

Since $\gcd(n,a)=1$ we have 
$\gcd(n,1-a)=\gcd(n,-a(1-a^{-1}))= \gcd(n,1-a^{-1})$
so that $n_a=n_{a^{-1}}$ as required. 
\end{proof}

\section{$\mathbb{Z}$-automorphisms and Computations}

Let $X$ be a finite Alexander quandle and let $X_A$ denote $X$ regarded as
an Abelian group, called the \textit{underlying Abelian group} of $X$. The 
map $\phi:X_A\to X_A$ defined by $\phi(x)=tx$ is
a homomorphism of $\mathbb{Z}$-modules. Since $t^{-1}\in \Lambda$, the
map $\psi:X_A\to X_A$ defined by $\psi(x)=t^{-1}x$ is a two-sided inverse
for $\phi$ as $\psi(\phi(x))=t^{-1}tx=x$ and $\phi(\psi(x))=tt^{-1}x=x$,
and $\phi$ is in fact a $\mathbb{Z}$-automorphism.

Conversely, if $A$ is a finite Abelian group and $\phi:A\to A$ is a
$\mathbb{Z}$-module automorphism, we can give $A$ the structure of an
Alexander quandle by defining $tx=\phi(x)$. This yields a general 
strategy for listing all finite Alexander quandles of a given size 
$n$: first, list all Abelian groups $A$ of order $n$; then, for each 
element of $\mathrm{Aut}_{\mathbb{Z}}(A)$ find $(1-t)A=\mathrm{Im}(1-\phi)$ 
and compare these as $\Lambda$-modules. In practice, for low order (i.e., 
$|A|\le 15$) Alexander quandles this procedure in its full generality is 
necessary only for one case, namely Alexander quandles with underlying 
Abelian group isomorphic to $\mathbb{Z}_4\oplus\mathbb{Z}_2$. We shall see
that Alexander quandles with $X_A\cong \mathbb{Z}_4\oplus\mathbb{Z}_2$ are 
isomorphic to linear Alexander quandles (in six cases) or to Alexander 
quandles with underlying group $(\mathbb{Z}_2)^3$ (in two cases).

We first obtain a few simplifying results:

\begin{lemma} \label{lem:cyc}
If the underlying Abelian group $X_A$ of $X$ is cyclic, then $X$ is linear.
\end{lemma}

\begin{proof}
Suppose $X_A=\mathbb{Z}_n$. Then for any $x\in \mathbb{Z}_n$ and any 
$\phi \in \mathrm{Aut}_{\mathbb{Z}}(\mathbb{Z}_n)$, we must
have $\phi(x)=\phi(x\cdot 1)=x\phi(1)$, so the action of $t$ agrees
with multiplication by $a=\phi(1)$ on $\mathbb{Z}_n$. Further, we must
have $\gcd(n,a)=1$ since $\phi$ is surjective. Hence $X$ is $\mathbb{Z}_n$ 
with $t$ acting by multiplication by $a$, that is, $X\cong\Lambda_n/(t-a)$. 
\end{proof}

\begin{remark}
Lemma \ref{lem:cyc} was also noted in \cite{O}.
\end{remark}

\begin{corollary} \label{cor:prime}
For any prime $p$, there are exactly $p-1$ distinct Alexander quandles
with $p$ elements, namely $\Lambda_p/(t-a)$ for $a=1,\dots,p-1$. Further,
every Alexander quandle of prime order is either trivial 
($\Lambda_p/(t-1)\cong T_p$, the trivial quandle of $p$ elements)
or connected.
\end{corollary}

\begin{proof}
If $p$ is prime, $n_a=\frac{n}{\mathrm{gcd}(p,1-a)}=1$ for each 
$a\in{1,\dots,p-1}$. Then by corollary \ref{cor:cond}, these are all 
distinct. By lemma \ref{lem:cyc}, every quandle of order $p$ 
is linear, so these are all of the Alexander quandles of order $p$.

Since $\gcd(p,1-a)=1$ for $a=2, \dots, p-1$, corollary \ref{cor:con2} 
gives us that $\Lambda_p/(t-a)$ is connected.
\end{proof}

\begin{corollary} \label{cor:prod}
Let  $n=p_1^{e_1}p_2^{e_2}\dots p_k^{e_k}$ be a product of powers of 
distinct primes. Then there are exactly $N_{p_1}N_{p_2}\dots N_{p_k}$ distinct 
Alexander quandles of order $n$, where $N_{p_i}$ is the number of distinct
Alexander quandles of order ${p_i}^{e_i}$. 
\end{corollary}

\begin{proof} 
Since any $\mathbb{Z}$-automorphism must respect order, any Alexander
quandle structure on a direct sum of Abelian groups 
$A_{p_1^{e_1}}\oplus\cdots\oplus A_{p_k^{e_k}}$ with order 
$p_1^{e_1},\dots, p_k^{e_k}$ must respect this direct sum structure. Hence 
we may obtain a complete list of Alexander quandles of order $n$ by 
listing all direct sums of Alexander quandles of orders 
$p_1^{e_1},\dots,p_k^{e_k}$.
\end{proof}

\begin{corollary} \label{cor:con4}
If the order of an Alexander quandle $n\equiv 2 (\mathrm{mod} \ 4)$, the
quandle is not connected.
\end{corollary}

\begin{proof}
If $n\equiv 2 (\mathrm{mod} \ 4)$, then the underlying Abelian group of
the quandle has a summand of $\mathbb{Z}_2$. Hence the quandle has a 
summand isomorphic to $\Lambda_2/(t+1)\cong T_2$, and therefore is not 
connected.
\end{proof}

\medskip

In light of corollary \ref{cor:prod}, to classify finite Alexander quandles 
it is sufficient to consider Alexander quandles of prime power order. 
Alexander quandles with prime order are cyclic as Abelian groups and hence 
are linear quandles, and so are classified by corollary \ref{cor:prime}. 
Alexander quandles with order a product of distinct primes are classified by 
corollary \ref{cor:prod}.

If the underlying Abelian group of $X$ is $(\mathbb{Z}_p)^n$, then
$X$ is not only a $\Lambda$-module but also a $\Lambda_p$-module, so
we may use the classification theorem for finitely generated modules over 
a PID. Thus any Alexander quandle $X$ with $X_A=(\mathbb{Z}_p)^n$ must be 
of the form $\Lambda_p/(h_1)\oplus\dots\oplus\Lambda_p/(h_k)$ with 
$h_1|h_2|\dots |h_k$, $h_i\in \Lambda_p$ and $\sum \mathrm{deg}(h_i)=n$. 
We may further assume without loss of generality that each 
$h_i\in \mathbb{Z}_p[t]$, is monic, and has nonzero constant term. 

\begin{prop} \label{prop:mconn}
An Alexander quandle $M=\Lambda_p/(h)$, $p$ a prime, is connected 
iff $(1-t)\not| h$.
\end{prop}

\begin{proof}
Since $M$ is finite, $(1-t)M=M$ iff $(1-t):M\to M$ is bijective. If 
$(1-t)|h$ then $h=(1-t)g$ for some nonzero $g\in M$, and hence 
$\mathrm{ker}(1-t)\ne \{0\} $, so $(1-t)$ fails to be injective.

Conversely, $(1-t)$ is prime in $\Lambda$, so $(1-t)$ coprime to $h$ implies 
that every $l\in \Lambda$ may be written as $a(1-t)+bh$ for some 
$a,b\in \Lambda$. Hence every $m\in M$ is $a(1-t)$ for some $a\in M$.
\end{proof}

\begin{prop}
The Alexander quandle 
$\Lambda_{p^n}/\left(t^n+\sum_{i=0}^{n-1}a_it^i\right)$ is 
connected iff $\sum_{i=0}^{n-1}a_i=-1$.
\end{prop}

\begin{proof}
By \ref{prop:mconn}, $\Lambda_{p^n}/(t^n+\sum_{i=0}^{n-1}a_it^i)$ is connected
iff $(t-1)|t^n+\sum_{i=0}^{n-1}a_it^i$. That is, 
$\Lambda_{p^n}/(t^n+\sum_{i=0}^{n-1}a_it^i)$ is connected iff there are 
$b_i \in \Lambda_p$, $0\le i\le n-2$ such that
\[
(t-1)\left(t^{n-1}+\sum_{i=0}^{n-2}b_it^i\right)=t^n+\sum_{i=0}^{n-1}a_it^i.
\]
Comparing coefficients, we must have that $a_{n-1}+b_{n-2}=-1$, 
$b_i=a_i+b_{i-1}$ for all $1\le i\le n-2$, and $b_0=a_0$. Then
$\sum_{i=0}^{n-1}a_i=-1$. Conversely, if $\sum_{i=0}^{n-1}a_i=-1$, define
$b_0=a_0$,  $b_i=a_i+b_{i-1}$ for all $1\le i\le n-2$, and 
$a_{n-1}+b_{n-2}=-1$.
\end{proof}

\begin{prop} \label{prop:p2}
There are $2p^2-3p-1$ connected Alexander quandles of order $p^2$ where
$p$ is prime.\footnote{This agrees with the result of Gra\~na in \cite{G}.}
\end{prop}

\begin{proof}
Every Alexander quandle of order $p^2$ has underlying Abelian group $\mathbb{Z}_{p^2}$ or
$\mathbb{Z}_p\oplus \mathbb{Z}_p$. A linear quandle $\Lambda_{p^2}/(t-a)$ 
of order $p^2$ is connected iff $\gcd(1-a,p)=1$, and there are $p(p-2)$ such 
quandles.

An Alexander quandle $M$ with underlying Abelian group 
$\mathbb{Z}_p\oplus \mathbb{Z}_p$ is a module over the PID $\Lambda_p$, 
so we have either $M\equiv \Lambda_p/(t-a)\oplus\Lambda_p/(t-a)$ or
$M\equiv\Lambda_p/(t^2+at+b)$ where $b\ne 0$. There are $p-2$ connected 
quandles of the first type and $(p-1)^2$ of the second type, so in total there 
are $2p^2-3p-1$ connected Alexander quandles of order $p^2$.
\end{proof}

For arbitrary values of $n$ and $p$ we may classify Alexander quandles 
with underlying Abelian group $(\mathbb{Z}_p)^n$ listing all possible 
$\Lambda$-modules with underlying group $(\mathbb{Z}_p)^n$ and comparing 
the submodules $\mathrm{Im}(1-t)$. 

Results of applying this procedure to Alexander quandles with underlying
Abelian group $(\mathbb{Z}_2)^2$, $(\mathbb{Z}_2)^3$ and $(\mathbb{Z}_2)^2$
are collected in Table 1. As we expect, these results agree with proposition
\ref{prop:p2}.

\begin{table}
\begin{center}
\begin{tabular}{|c|l|l|} \hline
$X_A$            & Module & $\mathrm{Im}(1-t)$ \\ \hline
 & $(\Lambda_2/(t+1))^2$ & 0 \\
$(\mathbb{Z}_2)^2$ & $\Lambda_2/(t^2+1)$   & $\Lambda_2/(t+1)$ \\
 & $\Lambda_2/(t^2+t+1)$ & $\Lambda_2/(t^2+t+1)$ \\ \hline
 & $(\Lambda_2/(t+1))^3$ & 0 \\
 & $\Lambda_2/(t+1) \oplus \Lambda_2/(t^2+1)$ & $\Lambda_2/(t+1)$ \\
$(\mathbb{Z}_2)^3$ & $\Lambda_2/(t^3+1)$ & $\Lambda_2/(t^2+t+1)$ \\
 & $\Lambda_2/(t^3+t+1)$ & $\Lambda_2/(t^3+t+1)$ \\
 & $\Lambda_2/(t^3+t^2+1)$ & $\Lambda_2/(t^3+t^2+1)$ \\
 & $\Lambda_2/(t^3+t^2+t+1)$ & $\Lambda_2/(t^2+1)$ \\ \hline
 & $(\Lambda_3/(t+2))^2$ & 0 \\
 & $(\Lambda_3/(t+1))^2$ & $(\Lambda_3/(t+1))^2$ \\
 & $\Lambda_3/(t^2+2)$ & $\Lambda_3/(t+1)$ \\
$(\mathbb{Z}_3)^2$ & $\Lambda_3/(t^2+1)$ & $\Lambda_3/(t^2+1)$ \\
 & $\Lambda_3/(t^2+2t+2)$ & $\Lambda_3/(t^2+2t+2)$ \\
 & $\Lambda_3/(t^2+2t+1)$ & $\Lambda_3/(t^2+2t+1)$ \\ 
 & $\Lambda_3/(t^2+t+2)$ & $\Lambda_3/(t^2+t+2)$ \\ 
 & $\Lambda_3/(t^2+t+1)$ & $\Lambda_3/(t+2)$ \\ \hline
\end{tabular}
\end{center}
\caption{Computations of $\mathrm{Im}(1-t)$ for $(\mathbb{Z}_2)^2$, 
$(\mathbb{Z}_2)^3$ and $(\mathbb{Z}_3)^2$.}
\end{table}

Note that by theorem \ref{thm:main} and corollary \ref{cor:cond}, the 
results in table 1 show that $\Lambda_2/(t^2+1)\cong\Lambda_4/(t-3)$ and 
$(\Lambda_2/(t+1))^2\cong \Lambda_4/(t-1)\cong T_4$, the trivial quandle 
of order 4, while $\Lambda_2/(t^2+t+1)$ is the only connected Alexander 
quandle of order 4.

Alexander quandles with underlying Abelian group $(\mathbb{Z}_2)^3$
include $\Lambda_2/(t+1)\oplus\Lambda_2/(t^2+1) \cong \Lambda_8/(t-5)$ 
and $(\Lambda_2/(t+1))^3\cong T_8$. Also, theorem \ref{thm:main} yields
an isomorphism $\Lambda_8/(t-3) \cong \Lambda_8/(t-7)$; otherwise, the 
order eight quandles listed are all distinct. Of these, only
$\Lambda_2/(t^3+t^2+1)$ and $\Lambda_2/(t^3+t+1)$ are connected. Note that
none of the linear Alexander quandles of order eight are connected.

Among Alexander quandles with Abelian group $(\mathbb{Z}_3)^2$, 
we have $\Lambda_9/(t-4) \cong \Lambda_9/(t-7)\cong \Lambda_9/(t^2+t+1)$ 
(the first isomorphism was noted in \cite{CJKS} and 
the second also follows from proposition 4.1 of \cite{LN}); otherwise, the 
linear quandles of order nine and the quandles listed in table 1 are all 
distinct. Note that five of the eight listed quandles of order nine are 
connected; of the linear quandles of order nine, $\Lambda_9/(t-2)$, 
$\Lambda_9/(t-5)$ and $\Lambda_9/(t-8)$ are connected.

To count distinct Alexander quandles whose underlying Abelian group is
neither cyclic nor a direct sum of $n$ copies of $\mathbb{Z}_p$, the 
following observation is useful.

\begin{lemma} \label{lem:ub}
The number of conjugacy classes in $\mathrm{Aut}_{\mathbb{Z}}(X_A)$ is an
upper bound on the number of distinct Alexander quandles $X$ with underlying
Abelian group $X_A$.
\end{lemma}

\begin{proof}
Let $\phi_1, \phi_2\in \mathrm{Aut}_{\mathbb{Z}}X_A$. Then if 
$t_1=\phi_1(1)$ and $t_2=\phi_2(1)$, we have $\phi_2^{-1}\phi_1\phi_2$ 
acting by multiplication by $t_2^{-1}t_1t_2=t_1$ since multiplication in 
$\Lambda$ is commutative. Thus any two conjugate automorphisms define the
same Alexander quandle structure.
\end{proof}

To complete the classification of Alexander quandles with up to fifteen
elements, we now only need to consider the case 
$X_A=\mathbb{Z}_4\oplus\mathbb{Z}_2$.

\begin{prop}
There are three distinct Alexander quandle structures definable on the Abelian group
$\mathbb{Z}_4\oplus\mathbb{Z}_2$, given by $\mathbb{Z}$-automorphisms
$\phi_1=\mathrm{id}$, $\phi_2\left((1,0)\right)=(1,1)$, 
$\phi_2\left((0,1)\right)=(0,1)$, $\phi_3\left((1,0)\right)=(1,1)$ and
$\phi_3\left((0,1)\right)=(2,1)$.
Further, these quandles are isomorphic to previously listed quandles, 
namely $(\mathbb{Z}_4\oplus\mathbb{Z}_2,\phi_1) \cong T_8$, 
$(\mathbb{Z}_4\oplus\mathbb{Z}_2,\phi_2) \cong 
\Lambda_2/(t+1)\oplus\Lambda_2/(t^2+1)$, and 
$(\mathbb{Z}_4\oplus\mathbb{Z}_2,\phi_3) \cong \Lambda_2/(t^3+t^2+t+1)$.
\end{prop}

\begin{proof}
Direct calculation shows that 
$\mathrm{Aut}_{\mathbb{Z}}(\mathbb{Z}_4\oplus\mathbb{Z}_2) \cong D_8$, the
dihedral group of order eight, so by lemma \ref{lem:ub} there are at 
most five Alexander quandle structures on $\mathbb{Z}_4\oplus\mathbb{Z}_2$. 
Of the eight $\mathbb{Z}$-automorphisms of $\mathbb{Z}_4\oplus\mathbb{Z}_2$, 
one is the identity, yielding the trivial quandle structure; five have 
$\mathrm{Im}(1-t)\cong \Lambda_2/(t+1)$ (including $\phi_2$) and hence 
yield quandles isomorphic to $\Lambda_2/(t+1)\oplus\Lambda_2/(t^2+1)$,
and two have $\mathrm{Im}(1-t)\cong \Lambda_2/(t^2+1)$ (including 
$\phi_3$), yielding quandles isomorphic to $\Lambda_2/(t^3+t^2+t+1)$. 
\end{proof}

We now have enough information to determine all Alexander quandles with
up to fifteen elements. In light of corollaries \ref{cor:prime} and 
\ref{cor:prod}, we list in table 2 only the numbers of distinct and
connected Alexander quandles of each order.

\begin{table}
\begin{center}
\begin{tabular}{|r|l|l|} \hline
 & \# of Alexander & \# \\ 
 $n$ & quandles & connected \\ \hline
 2 & 1  & 0 \\
 3 & 2  & 1 \\
 4 & 3  & 1 \\
 5 & 4  & 3 \\
 6 & 2  & 0 \\
 7 & 6  & 5 \\
 8 & 7  & 2 \\
 9 & 11 & 8 \\
10 & 4  & 0 \\
11 & 10 & 9 \\
12 & 6  & 1 \\
13 & 12 & 11 \\
14 & 6  & 0 \\
15 & 8  & 3 \\ \hline
\end{tabular}
\end{center}
\caption{The number of Alexander quandles and connected Alexander quandles 
of size $n\le 15$.}
\end{table}

\bibliographystyle{amsplain}

\end{document}